\documentclass{agtart_a}
\pdfoutput=1
\usepackage{pinlabel}

\title{Labeled binary planar trees and quasi-Lie algebras}

\author{Jerome Levine}
\givenname{Jerome}
\surname{Levine}
\address{Department of Mathematics\\
Brandeis University\\\newline
Waltham MA 02454-9110\\
USA\\\newline
\\\newline
{\rm Correspondence to:\qua K Orr\\\newline}
Department of Mathematics\\Indiana University\\\newline
831 East 3rd St\\Bloomington IN 47405-7106\\USA
}
\email{korr@indiana.edu}
\urladdr{http://php.indiana.edu/~korr/}

\volumenumber{6}
\issuenumber{}
\publicationyear{2006}
\papernumber{35}
\startpage{935}
\endpage{948}

\doi{}
\MR{}
\Zbl{}


\received{23 December 2005}
\revised{}
\accepted{2 May 2006}
\published{9 August 2006}
\publishedonline{9 August 2006}
\proposed{}
\seconded{}
\corresponding{}
\editor{}
\version{}

\dedicatory{This paper, Jerome Levine's fourth contribution to Algebraic
and Geometric Topology, is published posthumously, following the author's
untimely death in April 2006.  The editors are very grateful to Kent
Orr for preparing and proofreading the final version.}

\arxivreference{math.GT/0504278}

\makeatletter
\def\cnewtheorem#1[#2]#3{\newtheorem{#1}{#3}[section]
\expandafter\let\csname c@#1\endcsname\c@proposition}
\makeatother

\let\xysavmatrix\xymatrix
\def\xymatrix{\disablesubscriptcorrection\xysavmatrix}
\AtBeginDocument{\let\bar\wbar\let\tilde\wtilde\let\hat\what}

\def\draft{n}
\def\lbl#1{\label{#1}\printname{#1}}
\theoremstyle{plain}

\cnewtheorem{theorem}[proposition]{Theorem}
\cnewtheorem{lemma}[proposition]{Lemma}
\cnewtheorem{corollary}[proposition]{Corollary}
\cnewtheorem{claim}[proposition]{Claim}

\theoremstyle{definition}
\cnewtheorem{definition}[proposition]{Definition}
\cnewtheorem{problem}[proposition]{Problem}

\cnewtheorem{example}[proposition]{Example}
\cnewtheorem{exercise}[proposition]{Exercise}
\cnewtheorem{hint}[proposition]{Hint}
\cnewtheorem{remark}[proposition]{Remark}
\cnewtheorem{remarks}[proposition]{Remarks}

\renewcommand{\ni}{\noindent}

\renewcommand{\ker}{\operatorname{Ker}}
\newcommand{\im}{\operatorname{Im}}

\def\printname#1{
	\if\draft y
		\smash{\makebox[0pt]{\hspace{-0.5in}
			\raisebox{8pt}{\tt\tiny #1}}}
	\fi
}

\def\SS{\Sigma}

\def\co{\colon\thinspace}

\def\a{\alpha}
\def\b{\beta}
\def\d{\delta}

\def\o{\omega}

\def\g{\gamma}
\def\l{\lambda}

\def\Z{\mathbb Z}

\def\D{\mathsf D}
\def\H{\mathcal H}
\def\C{\mathcal C}
\def\Dr{\mathcal D}

\def\B{\mathcal B}

\def\L{{\mathsf L}}

\def\Q{\mathbb Q}

\def\sub{\subseteq}

\def\T{\mathcal T}

\def\ti{\tilde}

\def\sf#1{\mathcal S^{\text{fr}}_{#1}} 
\def\At#1{\mathcal A^t_{#1}(H)}
\def\DH#1{\mathsf D_{#1}(H)}
\def\LH#1{\mathsf L_{#1}(H)}

\begin{document}

\begin{asciiabstract}
We study the natural map eta between a group of binary planar
trees whose leaves are labeled by elements of a free abelian group H
and a certain group D(H) derived from the free Lie algebra
over H. Both of these groups arise in several different topological
contexts. The map eta is known to be an isomorphism over Q, but not
over Z. We determine its cokernel and attack the conjecture
that it is injective.
\end{asciiabstract}

\begin{htmlabstract}
We study the natural map &eta; between a group of binary planar
trees whose leaves are labeled by elements of a free abelian group H
and a certain group D(H) derived from the free Lie algebra
over H. Both of these groups arise in several different topological
contexts. &eta; is known to be an isomorphism over <b>Q</b>, but not
over <b>Z</b>. We determine its cokernel and attack the conjecture
that it is injective.
\end{htmlabstract}

\begin{abstract}
We study the natural map $\eta$ between a group of binary planar
trees whose leaves are labeled by elements of a free abelian group $H$
and a certain group $\mathsf{D}(H)$ derived from the free Lie algebra
over $H$. Both of these groups arise in several different topological
contexts. $\eta$ is known to be an isomorphism over $\mathbb{Q}$, but not
over $\mathbb{Z}$. We determine its cokernel and attack the conjecture
that it is injective.
\end{abstract}

\maketitle

\section{Introduction}
Let $H$ be a finitely-generated free abelian group and $\L (H)$ the graded
free Lie algebra on $H$. There is a natural homomorphism $H\otimes \L
(H)\to\L(H)$ defined by bracketing, whose kernel is denoted $\D (H)$. If
$H$ supports a non-singular symplectic form, eg $H=H_1 (\SS )$, where
$\SS$ is a closed orientable surface, with symplectic basis $\{ x_i
,y_i\}$, then $\D (H)$ is, in fact, a Lie algebra. It can be identified
with the Lie subalgebra of $\Dr (\L (H))$ (the graded Lie algebra of
derivations of $\L (H)$) consisting of those derivations  which vanish
on the element $\sum_i [x_i ,y_i ]\in\L_2 (H)$.

$\D (H)$ has arisen in several different topological contexts. For
example, it was probably first observed by Orr \cite{O} (but see also
Habegger and Lin \cite{HL}) that it is very natural to regard the Milnor invariants of
a link $L$, or, more precisely, string link, as elements of $\D (H)$,
where $H=H_1 (S^3 -L)$. If $\SS$ is a compact orientable surface with
one boundary component and $H=H_1 (\SS )$ then $\D (H)$ contains, as a
Lie subalgebra, the associated graded Lie algebra of the {\em relative
weight} filtration, defined by D. Johnson, of the mapping class group
of $\SS$ -- see Johnson \cite{J} and Morita \cite{M}. Similarly, if we
consider the homology  concordance group of homology cylinders over a
surface -- see Garoufalidis--Levine \cite{GL} and Levine \cite{L}
-- there is also a relative weight filtration and, in this case, the
associated graded Lie algebra is actually isomorphic to $\D (H)$. $\D
(H)$ appears in Kontsevich's work \cite{K} on graph complexes and his
computation of the cohomology of the group of outer automorphisms of a
free group.

Consider now the abelian group $\At{}$ generated by unitrivalent trees,
with cyclic orientations of its trivalent vertices and univalent vertices
labelled by elements of $H$, subject to the anti-symmetry and IHX
relations and linearity of the labels. $\At{}$ appears as the indexing
of the so-called {\em tree-level} of the Kontsevich integral of a link
or string link. See Habegger and Masbaum \cite{HM} where this is related
to the Milnor invariants via a natural map $\eta\co\At{}\to\DH{}$. It
is proved there that, rationally, the Milnor invariants of a string link
determine the tree-level of its Kontsevich integral. This corresponds to
the fact that the map $\eta\otimes\Q \co \At{}\otimes\Q\to\DH{}\otimes\Q$
is an isomorphism, which is proved by Habegger and Pitsch \cite{HP} (see
also Garoufalidis--Levine \cite{GL} and Levine \cite{L1}). $\At{}$
appears in Habiro \cite{H} and the study by Garoufalidis, Goussarov and
Polyak \cite{GGP} of claspers and finite-type invariants of 3-manifolds,
and subsequently in Levine \cite{L}, mapping onto the associated
graded groups of a filtration of the concordance group of homology
cylinders defined using {\em claspers}. In this context the map $\eta
\co \At{}\to\DH{}$ reflects the relation between the clasper filtration
and the usual relative weight filtration. Most recently the group $\At{}$
appears in the work of Schneiderman--Teichner \cite{ST}, where it encodes
the obstruction to removing intersection and self-intersection points
of immersed connected surfaces in a simply-connected 4-manifold (here
$H$ is a free abelian group of rank equal to the number of surfaces)
via a tower of Whitney disks. In the special case where the surfaces are
disks in the 4-ball bounded by a link in $S^3$, the obstruction element
in $\At{}$ maps to the element in $\DH{}$ corresponding to the Milnor
invariants of the link.

In these various situations the study of the homomorphism $\At{}\to\DH{}$
is closely related to the question of whether there are invariants
in $\At{}$ which give more information than the analogous, perhaps
more easily defined, invariants in $\DH{}$. For example the work of
Schneiderman--Teichner may uncover new invariants of link concordance
beyond the Milnor invariants.

Since $\eta\otimes\Q$ is an isomorphism, the kernel and cokernel of
$\eta$ are finite. In \cite{L1} some progress was made toward determining
them. Toward this end we introduced the notion of a {\em quasi-Lie
algebra} and studied the structure of a free quasi-Lie algebra. In the
present work we extend these results. In particular we determine the
precise structure of a free quasi-Lie algebra (adapting an argument of
Marshall Hall \cite{Ha}), determine the cokernel of $\eta$ precisely
and show that a ``good part'' of the kernel of $\eta$ is trivial. It
remains a reasonable conjecture that $\eta$ is injective.

It has recently come to our attention that some of the results of this
note (and of \cite{L1}) have been independently obtained by K Habiro --
in particular \fullref{th.at}, \fullref{th.etao} in \fullref{seceta}
(and \fullref{cor.prime}) by a very similar method, as well as
\fullref{th.eta2} by a different method.

The author was partially supported by an NSF grant.

\section{Statement of results}

We will use the precise definition of the groups $\{\DH{n}\}$ and
$\{\At{n}\}$ given in \cite{L,L1}. In particular $\DH{n}$ is the
kernel of the bracket map $H\otimes\L_{n+1}(H)\to\LH{n+2}$ and $\At{n}$
is a quotient of the free abelian group generated by unitrivalent trees
with $n$ trivalent vertices, each of which is given an orientation --
ie a cyclic ordering of its incident edges -- and whose univalent
vertices are labeled by elements of $H$. The ``relations'' which are
divided out are:
\begin{itemize}
\item  anti-symmetry: $T+T'=0$, where $T'$ is identical with $T$ except that one trivalent vertex is given the opposite orientation.
\item IHX: $T_1 -T_2 +T_3 =0$, where the $T_i$ are identical except in the neighborhood of two adjacent trivalent vertices, which look as follows:

\medskip
\centerline{
\labellist\small
\pinlabel {$T_1$} at 40 20
\pinlabel {$T_2$} at 160 20
\pinlabel {$T_3$} at 275 20
\endlabellist
\includegraphics[width=3in]{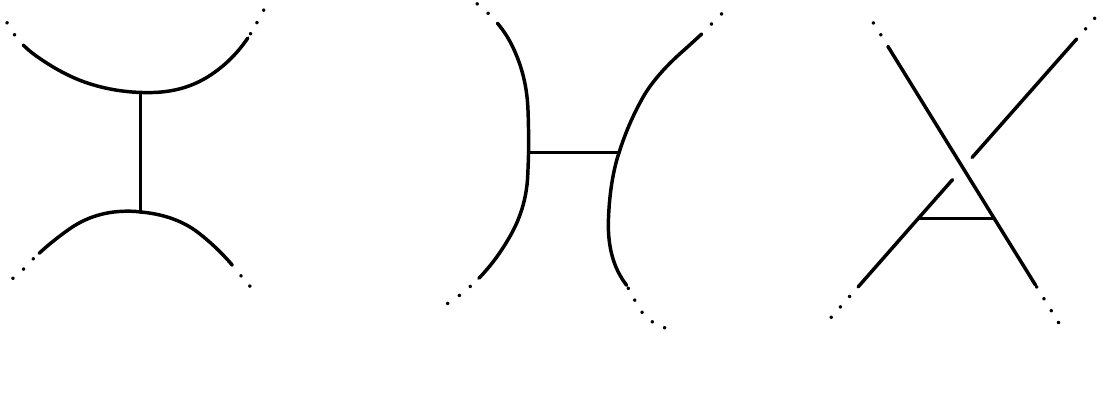}
}

The orientations of the trivalent vertices are counterclockwise.
\item linearity: $T=T_1 +T_2$ where $ T_1 ,T_2 , T $ are identical except that one of the univalent vertices has labels $a_1 ,a_2 ,a_1 +a_2$, respectively, for some $a_i\in H$.
\end{itemize}

The graphical representation of free Lie algebras over $\Q$, for example,
is well-known. Extending this to Lie algebras over $\Z$ requires some
extra considerations. In \cite{L1} we introduce the notion of a quasi-Lie
algebra, in which the relation $[\a ,\a ]=0$ is replaced by the slightly
weaker relation $[\a ,\b ]+[\b ,\a ]=0$. Then the free quasi-Lie algebra
$\L' (H)$ over a free abelian group $H$ is isomorphic to a Lie algebra
of trees similar to the definition of $\At{}$ except that one univalent
vertex (the {\em root}) is not labeled.

There is an obvious epimorphism $\g_n \co \L_n '(H)\to\L_n
(H)$. We can also define a bracketing homomorphism $\b_n ' \co
H\otimes\L_{n+1}'(H)\to\L_{n+2}'(H)$ by $\b_n '(h\otimes\l )=[h,\l ]$
and then define $\DH{n}=\ker\b_n '$. A map $\eta'_n \co \At{n}\to
H\otimes\L_{n+1}'(H)$ is defined by
\begin{equation}\lbl{eq.eta}
\eta'_n (T)=\sum_i h_i\otimes [T_i ]
\end{equation}
where the summation is over all univalent vertices of $T$. For each univalent vertex $h_i$ is its label and $T_i$ is the rooted tree obtained from $T$ by making that vertex the root -- $[T_i ]$ is the corresponding element of $\L'_{n+1}(H)$.

In \cite{L1} the following theorem is proved.
\begin{theorem}\lbl{th.at}
The sequence
$$
\At{n}\xrightarrow{\eta'_n}H\otimes\L'_{n+1}(H)\xrightarrow{\b'_n}\L'_{n+2}(H)\to 0
$$
is exact.
\end{theorem}
Therefore $\eta'_n$ defines an epimorphism $\At{n}\to\D'_n (H)$.

In order to completely understand the original map $\At{n}\to\DH{n}$ we need to resolve the following two problems:
\begin{enumerate}
\item Determine the map $\g_n \co \D'_n (H)\to\DH{n}$, ie determine its kernel and cokernel.
\item Determine the kernel of $\eta'_n$.
\end{enumerate}
Problem (1) essentially reduces to determining $\L'_n (H)$, since $\L_n (H)$ is well-understood.

In \cite{L1} it is shown that $\g_n$ is an isomorphism if $n$ is odd and,
for $n$ even there is an exact sequence
$$
\L_k (H)\otimes\Z/2\xrightarrow{\theta_k}\L'_{2k}(H)\xrightarrow{\g_{2k}}\L_{2k}(H)\to 0
$$
where $\theta_k$ is defined by $\theta_k (\a )=[\a ,\a ]$ (where
$\a\in\L_k (H)$ is lifted into $\L'_k (H)$). It was conjectured in
\cite{L1} that $\theta_k$ is injective. Our first result proves this
conjecture.
\begin{theorem}\lbl{th.fla}
The sequence
$$
0\to\L_k (H)\otimes\Z/2\xrightarrow{ \theta_k}\L'_{2k}(H)\xrightarrow{\g_{2k}}\L_{2k}(H)\to 0
$$
is exact (and therefore split exact).
\end{theorem}
As a consequence of this theorem, we have the following relations between $\D'_n (H)$ and $\DH{n}$.
\begin{corollary}\lbl{cor.dd}
There exist exact sequences:
$$0\to\D'_{2k}(H)\to\DH{2k}\to\L_{k+1}(H)\otimes\Z/2\to 0 $$
$$0\to H\otimes\L_k (H)\otimes\Z/2\to\D'_{2k-1}(H)\to\DH{2k-1}\to 0$$
\end{corollary}
These exact sequences are derived in \cite{L1}.

\begin{remark} We can describe the elements of $\DH{2k}$ which do not come from $\D'_{2k}(H)$ in the following graphical manner. Let $\a\in\L_{k+1}(H)$ be represented by a labeled rooted tree $T$ of degree $k$. Let $T'$ be another copy of $T$ and let $T\odot T'$ be the labeled tree (representing an element of $\At{2k}$) obtained by welding the roots of $T$ and $T'$ together. For each labeled univalent vertex $v_i$ of $T$ let $T_i\odot T'$ be the rooted labeled tree obtained from $T\odot T'$ by making $v_i$ the root. If $h_i$ is the label of $v_i$ in $T$ then consider the element $\sum_i h_i\otimes (T_i\odot T')\in H\otimes\L'_{2k+1}(H)$. This does not lie in $\D'_{2k}(H)$ but its projection into $H\otimes\L_{2k+1}(H)$ does lie in $\DH{2k}$ and maps to $\a\otimes 1\in\L_{k+1}(H)\otimes\Z/2$.

\end{remark}
We now turn to Problem (2). In \cite{L1} the following is proved.
\begin{theorem}\lbl{th.eta}
$\eta'_n$ is a split surjection. $\ker\eta'_n$ is the torsion subgroup of $\At{n}$ if $n$ is even, and is the odd torsion subgroup of $\At{n}$ if $n$ is odd.

In both cases $(n+2)\ker\eta'_n =0$.
\end{theorem}
One immediate consequence is the known result that $\At{n}\otimes\Q\cong\DH{n}\otimes\Q$.

We will improve on this result by constructing a splitting of $\At{}$ and $\D'(H)$ such that $\eta'$ preserves components, and give a better estimate on the order of the kernel of each factor. In particular we will show:
\begin{corollary}\lbl{cor.prime} If $n+2$ is a prime power $p^k$, then $p^{k-1}\ker\eta'_n =0$. For example,
if $n+2$ is prime, then $\eta'_n$ is an isomorphism.
\end{corollary}

Finally, by a direct computation of ranks we will show:
\begin{theorem}\lbl{th.eta2}
$\eta'_2$ is an isomorphism.
\end{theorem}
Since it is obvious that $\eta'_1$ is an isomorphism, the first unsettled case is $n=4$.

\section[Structure of the quasi-Lie algebra: Proof of Theorem~\ref{th.fla}]
{Structure of the quasi-Lie algebra: Proof of \fullref{th.fla}}

Choose a basis $\{ a_1 ,\ldots ,a_m\}$ of $H$. Let $\C_n$ denote the set
of formal commutators of degree $n$ in the $a_i$ and $\C =\cup_n\C_n$.
Recall the definition of a {\em Hall basis} (see, for example, Hall
\cite{Ha}). Choose a linear ordering of the elements of $\C$ satisfying
only that if $d(x)> d(y)$ (where $d$ denotes degree), then $x>y$. Let $\H$
be the subset of $\C$ defined recursively by the following properties:
\begin{enumerate}
\item Each $a_i\in\H$
\item If $u,v\in\C$, then $[u,v]\in\H$ if and only if:
\begin{enumerate}
\item $u, v\in\H$
\item $u>v$
\item If $u=[x,y]$ (and so $x,y\in\H$ and $x>y$), then $v\ge y$.
\end{enumerate}

\end{enumerate}
Note that $\H$ depends on the choice of ordering.

It is a well-known result (see, for example, \cite{Ha}) that any Hall
basis is a basis of the free Lie algebra $\LH{})$.

Let $\ti\H$ denote the subset of $\C$ consisting of all elements of the form $[h,h]$ for some $h\in\H$. It is clear that \fullref{th.fla} will follow from:
\begin{lemma}\lbl{lem.fla}
$\L'(H)\otimes\Z/2$ has, as basis, $\H'=\H\cup\ti\H$.
\end{lemma}
\begin{proof}
We will follow closely the proof in \cite{Ha}, making a few necessary
modifications to apply to our situation.

Let $V_n$ be the $\Z/2$--vector space with basis $\C_n$, $V=\oplus_n V_n$,
and $W_n$ the $\Z/2$--vector space with basis $\H'_n$ and $W=\oplus_n W_n$. There are obvious maps:
$$W_n \sub V_n\to\L'_n (H)$$

We will define a retraction $r \co V_n\to W_n$ recursively on $n$, satisfying
\begin{enumerate}
\item If $h\in\H'$ then $r(h)=h$.
\item For any $c\in\C$,  $r(c)=c$ in $\L'(H)$.
\item For any $c_1 ,c_2\in\C$, $r[c_1 ,c_2 ]=r[r(c_1 ),r(c_2 )]$.
\end{enumerate}
For $n=1$ we define $r(a_i )=a_i$.

Now suppose $r$ is defined on $V_k$ for all $k<n$ satisfying (1)--(3). We will define a sequence of additive moves $V_n\to V_n$ which will define $r$ when it stops.

\ni \textbf{Step 1}\qua If $c=[c_1 ,c_2 ]$, then $c\to [r(c_1 ), r(c_2 )]$.

Now apply Step 2 to each term of the sum.

\ni \textbf{Step 2}\qua If $c=[h_1 ,h_2 ]$, where $h_1 ,h_2\in\H'$, then
$$c\to\begin{cases}
0 & \text{if $h_1$ or $h_2$ belongs to $\ti\H$ (case 1)}\\
c & \text{if $h_1 ,h_2\in\H$ and $h_1 \ge h_2$ (case 2)}\\
[h_2 ,h_1 ] &\text{if $h_1 ,h_2\in\H$ and $h_1 <h_2$ (case 3)}

\end{cases}$$

In case 1 stop. In case 2 or 3 go on to Step 3.

\ni \textbf{Step 3}\qua If $c=[h_1 ,h_2 ]$, with $h_i\in\H$ and $h_1\ge h_2$, write $h_1 =[h_3 ,h_4 ]$ (note $h_3 >h_4$). Then
$$c\to\begin{cases}
c & \text{if $h_2\ge h_4$ (case 1)}\\
[[h_3 ,h_2 ], h_4 ]+[[h_2 ,h_4 ], h_3 ]& \text{if $h_2 <h_4$ (case 2)}
\end{cases} $$

In case 1 stop. In case 2, apply Step 1 to each of the terms in the sum.

It is clear that if $c\in\H'$, then the process will stop at Step 2 or 3
at $c$. In general we need to show that this process will stop after a
finite number of steps. It is clear then that properties (1)--(3) will be satisfied.

Define  a new relation among the elements of  $\C_n$, for $n\ge 2$. Let $c=[c_1 ,c_2 ]$ and $c'=[c'_1 ,c'_2 ]$. We will say 
$$c\succ c' \text{ if } \min (c_1 ,c_2 )>\min (c'_1 ,c'_2 ) 
$$

Now if $b$ is one of the terms in $[r(c_1 ) ,r(c_2 )]$ obtained after Step
1, then applying Steps 2 and 3 to $b$ will either stop, resulting in an
element of $\H'$ or lead us to case 2 of Step 3.
 Take $[c'_1 ,c'_2 ]$ to be either of the resulting terms in case 2 of
Step 3. When we then apply Step 1 we have $[r(c'_1 ), r(c'_2 )]\succ b$,
since $r(c'_2 )=c'_2 >h_2$ and $d(r(c_1 ))=d(c_1 )>d(h_2 )$. Thus
iterating the process results in a sum of terms each of which stabilizes
or leads to a sum of terms which are greater under the relation $\succ$.
Since there are only a finite number of elements in $\C_n$ the process
must eventually stop. In fact it must stop whenever the element $[h_1 ,h_2
]$ to which we apply Step 3 satisfies $d(h_2 )>\tfrac{n}{3}$ since this will force $d( h_4 )<\tfrac{n}{3}$ from which it follows that $h_2 >h_4$.

We now have defined a retraction $r \co V_n\to W_n$ satisfying properties
(1)--(3). To complete the proof we need to show that $r$ induces a map
$L'_n (H)\to W_n$. Notice that we can regard $V$ as the {\em free $\Z/2$--magma over $H$} and that $L' (H)\otimes\Z/2$ is the quotient of $V$ by the ideal $I$ generated by elements of the form 
$$\xi =[x,y]+[y,x] ]\text{\quad or\quad} [[x,y],z]+[[y,z],x]+[[z,x],y]$$
where $x,y,z\in\C$. So we need to show that $r(I)=0$.

Now $I$ is generated additively by formal brackets of elements of $V$, one of which is of the form $\xi$ above. By property (3) of $r$ it is only necessary then to show that $r(\xi )=0$. In fact, again by property (3), we may assume that $x,y,z\in\H'$.

\ni ${\xi =[x,y]+[y,x]}$:\qua
We may assume $x\ge y$. In the definition of $r$, we see that Step 2 will change $\xi$ either to $0$, if $x$ or $y$ belongs to $\ti\H$, or to $2[x,y]=0$ otherwise.

\ni $\xi =[[x,y],z]+[[y,z],x]+[[z,x],y]$:\qua
We may assume $x\ge y\ge z$ and proceed by induction on $d([[x,y],z])$.

\ni{\bf Case 1}\qua $x=y$:\qua
So $\xi =[[x,x],z]+[[x,z]+[z,x],x]$. Now $[x,x]$ either belongs to $\ti\H$ (if $x\in\H$) or $r([x,x])=0$ if $x\in\ti\H$. In either case $r([[x,x],z])=0$.

For the remaining terms note that $[x,z]+[z,x]\to 2[x,z]=0$.

\ni{\bf Case 2}\qua $x\in\ti\H$:\qua
Then $r([x,y])=r([[y,z],x])=r([z,x])=0$ by Step 2.

\ni{\bf Case 3}\qua $[x,y]\in\H$:\qua
In evaluating $r([[x,y],z])$ we proceed to Step 3 and apply case~2:
$$[[x,y],z]\to [[x,z],y]+[[z,y]x]$$
At this point $\xi$ has been reduced to $0$.

\ni{\bf Case 4}\qua $[x,y]\notin\H'$ and $x\in\H$:\qua
We proceed by a downward lexicographical induction. Assume that $r(\xi )=0$ when 
$$\xi = [[x',y'],z']+[[y',z'],x']+[[z',x'],y']$$
and $x'\ge y'\ge z'$ and either $z'>z$ or $z'=z$ and $y'>y$.

Write $x=[x_1 ,x_2 ]$. Since $x\in\H$, then $x_1 ,x_2\in\H$ and $x_1 >x_2$. Since $[x,y]\notin\H$ , then $x_2 >y$.

We therefore have:
\begin{align*}
r([[x,y],z])&=r([[[x_1 ,x_2 ],y],z])\\
&=r([([x_1 ,y],x_2 ]+[[x_2 ,y],x_1 ]),z])
\end{align*}
by Step 3, case 2, since $x_2 >y$, and Step 2, case 3, for the first term, since $x_1 >y$.

We can apply our downward induction to both terms, since $[x_1 ,y]>x_2 >y$ and both $[x_2 ,y]$ and $x_1$ are $>y$ (using Step 2 case 3, if $x_1 >[x_2 ,y]$) to obtain:
\begin{multline}\lbl{eq.hall}
r([[x,y],z])=r([[x_2 ,z],[x_1 ,y]]+[[[x_1 ,y],z],x_2 ]\\
+[[x_1 ,z],[x_2 ,y]]+[[[x_2 ,y],z],x_1 ])
\end{multline}

By our ongoing induction on $\deg ([[x,y],z])$ we have
\begin{align*}
r([[x_1 ,y],z])&=r([[x_1 ,z],y]+[[y,z],x_1 ])\\
r([[x_2 ,y],z])&=r([[x_2 ,z],y]+[[y,z],x_2 ])
\end{align*}
Substituting these equalities into equation \eqref{eq.hall} gives
\begin{multline}\lbl{eq.hall1}
r([[x,y],z])=r([[x_2 ,z],[x_1 ,y]]+[[[x_1 ,z],y], x_2 ]+[[[y,z],x_1 ],x_2 ]\\
+[[x_1 ,z],[x_2 ,y]]+[[[x_2 ,z],y],x_1 ]+[[[y,z],x_2 ],x_1 ])
\end{multline}
Now we write
\begin{align*}
r([[z,x],y])&=r([[z,[x_1 ,x_2 ]],y])=r([[[x_1 ,x_2 ],z],y])  \\
 &=r([[[x_1 ,z],x_2 ],y]+[[[x_2 ,z],x_1 ],y])
\end{align*}
using the Jacobi identity on elements of degree $<\deg ([[x,y],z])$.

We can now use our downward induction on each of the two terms on the right to get:
\begin{multline}\lbl{eq.hall2}
r([[z,x],y])=r([[x_1 ,z],[x_2 ,y]]+[[[x_1 ,z],y],x_2 ]+[[x_2 ,z],[x_1 ,y]]\\
+[[[x_2 ,z],y],x_1 ])
\end{multline}
We can now add equations \eqref{eq.hall1} and \eqref{eq.hall2}, cancelling out many of the terms, to get
\begin{multline}
r([[x,y],z]+[[y,z],x]+[[z,x],y])=r([[[y,z],x_1 ],x_2 ]\\
+[[[y,z],x_2 ],x_1 ] +[[y,z],[x_1 ,x_2 ]])
\end{multline}
If $[y,z]>x_2$, then our downward induction, applied to $[[[y,z],x_1 ],x_2]$, will tell us that $r([[[y,z],x_1 ],x_2 ]+[[[y,z],x_2 ],x_1 ]+[[y,z],[x_1 ,x_2 ]])=0$. (In case $x_1 <[y,z]$, we use Step 2.) If $[y,z]<x_2$ then we apply downward induction on $[[x_1 ,x_2 ],[y,z]]$.

Finally note that, since $x_1 >x_2 >y>z$, we conclude that 
$$
\deg z\le\tfrac14\deg([[x,y],z])=\tfrac{n}{4},
$$ 
which ensures  the beginning of the  induction.

This completes the proof of \fullref{lem.fla} and \fullref{th.fla}.
\end{proof}

\section{Study of $\eta'$}\label{seceta}

\subsection{A splitting of $\At{n}$}
We now consider the maps $\eta'_n \co \At{n}\to\D'_n (H)$. In \cite{L1}
it is proved that $$(n+2)\ker\eta'_n =0.$$ We will construct a splitting
of the various groups
$$
\At{}, \L'(H), H\otimes \L'(H)
$$
so that $\eta'$ will respect the summands of these splittings and then give better estimates of the order of $\ker\eta'$ on each summand. The splitting will depend on the choice of a basis $\B =\{\a_1 ,\ldots ,\a_d\}$ of $H$. We adopt a slightly different, but equivalent, view of $\At{}$ as generated by vertex-oriented unitrivalent trees with univalent vertices labeled by elements of $\B$, subject to the anti-symmetry and IHX relations (but now the linearity relation is not needed). Similarly $\L'(H)$ is generated by formal brackets in the elements of $\B$, subject to anti-symmetry and Jacobi relations. Thus $\L'(H)$ is graphically described as generated by vertex-oriented unitrivalent trees with univalent vertices labeled by elements of $\B$ and one unlabeled univalent vertex chosen as a ``root'', subject to anti-symmetry and IHX.

Let $\o =(n_1,\ldots ,n_d )$ be a sequence of non-negative integers. We will say that a labeled vertex-oriented unitrivalent tree has {\em signature $\o$} if exactly $n_i$ of the vertices are labeled by $\a_i$. A formal bracket has signature $\o$ if exactly $n_i$ of the entries in the bracket are $\a_i$. Notice that each anti-symmetry, IHX or Jacobi relation is defined by a sum of trees or brackets which all have the same signature.

We now define $\At{\o}$ to be the abelian group generated by labeled vertex-oriented unitrivalent trees of signature $\o$, subject to the anti-symmetry and IHX relations, and $\L'_{\o}(H)$ to be the group generated by brackets of signature $\o$, subject to anti-symmetry and Jacobi relations. It is clear that 
$$\At{n}=\oplus_{\o}\At{\o}\qquad \L'_n (H)=\oplus_{\o}\L'_{\o}(H)$$
where the sums range over all $\o$  with $\sum_i n_i =n+2$, for $\At{n}$ and with $\sum_i n_i =n$ for $\L'_n (H)$. Note that the Lie bracket defines a pairing
$$\L'_{\o}(H)\otimes\L'_{\o'}(H)\to\L'_{\o +\o'}(H)$$
where, if $\o =(n_1 ,\ldots ,n_d )$ and $\o'=(n'_1 ,\ldots ,n'_d )$, then $\o +\o' =(n_1 +n'_1 ,\ldots ,n_d +n'_d )$.

The map $\eta'_n \co \At{n}\to H\otimes\L'_{n+1}(H)$ is defined in equation \eqref{eq.eta}. If the tree $T$ in that formula has signature $\o=(n_1, \ldots ,n_d )$, then each term on the right side will be of the form $\a_j\otimes [T_i ]$, for some $j$, and where $[T_i ]\in\L' _{n+1}(H)$
has signature $\o_j =(n_1 ,\ldots ,n_{j}-1,\ldots ,n_d )$. Therefore we can write
\begin{equation*}
\eta'_n (T)=\sum_j \a_j\otimes\l_j \quad\text{ where }\l_j\in\L'_{\o_j}(H)
\end{equation*}
Now define $(H\otimes\L'(H))_{\o}=\sum_j \a_j\otimes\L'_{\o_j}(H)$. It is clear that 
$$
H\otimes\L'_{n+1}(H)=\oplus_{\o}(H\otimes\L' (H))_{\o},
$$
where $\o$ ranges over all $\o$
with $\sum_i n_i =n+2$ and $\eta'_n (\At{\o})\sub (H\otimes\L'(H))_{\o}$.

If $\o =(n_1 ,\ldots ,n_d )$ we define $\d (\o )=$greatest common divisor of $n_1 ,\ldots ,n_d$.
\begin{theorem}\lbl{th.etao}
$$\d (\o )\left(\ker\eta'_n\left|\At{\o}\right.\right)=0$$
\end{theorem}
\begin{proof}
In the proof of \fullref{th.at} in \cite{L1} we use a map $\rho_n
\co H\otimes \L'_{n+1}(H)\to\At{n}$ which sends any labeled tree
with a root to the same tree, forgetting which vertex is the root. The
observation that $\rho_n\circ\eta'_n=$multiplication by $n+2$ shows that
$(n+2)\ker\eta'_n =0$.

Now it is clear that $\rho_n (H\otimes\L'(H))_{\o}\sub\At{\o}$. But the restriction of $\rho_n$ to $(H\otimes\L'(H))_{\o}$ can be decomposed into a sum of maps $\rho_n^{\o ,i} \co (H\otimes\L'(H))_{\o}\to\At{\o}$ defined by
$$\rho_n^{\o ,i}(\a_j\otimes\l )=\begin{cases}
\rho_n (\a_j\otimes\l)&\text{ if } i=j\\
0&\text{ if } i\not=j
\end{cases}$$
It is clear that $\rho_n^{\o ,i}\circ\eta'_n\left|\At{\o}\right.$ is just multiplication by $n_i$, and so 
$$
n_i\ker\eta'_n\left|\At{\o}\right.=0
$$ 
for $i=1,\ldots ,d$.

This completes the proof.
\end{proof}
\subsection[Proof of Corollary~\ref{cor.prime}]{Proof of \fullref{cor.prime}}
We only need show that, for any $\o=(n_1 ,\ldots,n_d)$ with $\sum_i n_i =p^k$, $\d (\o )\left|p^{k-1}\right.$. Clearly $\d (\o )\left|p^k\right.$, so suppose $\d (\o )=p^k$. This can only happen if some $n_i =p^k$ and the remaining $n_j =0$, which means that every tree $T$ in the generating set of $\At{\o}$ has all its univalent vertices labeled by $\a_i$. Choose two univalent vertices which are each connected by an edge to the same trivalent vertex. Unless $n=1$ this trivalent vertex is connected by its third edge to another trivalent vertex. If we apply the IHX relation here, we see that $T=0$ in $\At{}$.

\medskip
\centerline{
\hair=1pt
\labellist\small
\pinlabel {$\alpha_i$} [br] at 27 150
\pinlabel {$\alpha_i$} [bl] at 60 150
\pinlabel {$\alpha_i$} [br] at 125 150
\pinlabel {$\alpha_i$} [bl] at 153 150
\pinlabel {$\alpha_i$} [br] at 225 150
\pinlabel {$\alpha_i$} [bl] at 260 150
\normalsize
\pinlabel {$-$} at 95 105
\pinlabel {$+$} at 200 105
\pinlabel {$= 0$} at 315 105
\endlabellist
\includegraphics[width=3in]{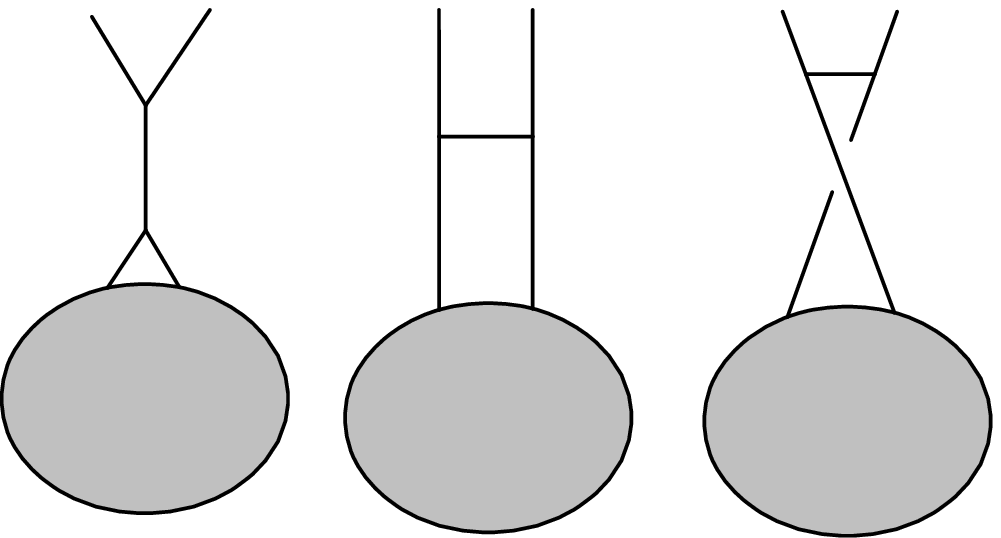}
}

In case $n=1$, we see, by anti-symmetry, that $2T=0$. But since $\d (\o )=3$ we also have $3T=0$.

\subsection[Proof of Theorem~\ref{th.eta2}]{Proof of \fullref{th.eta2}}
If follows from \fullref{th.at} and \fullref{cor.dd} that $\D'_2
(H)=\im\eta'_2$ is a free abelian group of rank $dd_3 -d_4$, where
$d_n$ denotes the rank of $\L_n (H)$. Therefore in order to prove that
$\eta'_2$ is injective it suffices to show that $\At{2}$ is generated by
$dd_3 -d_4$ elements.  Witt's formula (see, for example, Magnus, Karass
and Solitar \cite[Theorem 5.11]{MKS})
gives a general formula for $d_n$ -- in particular $d_3 =2\binom{d+1}{3}$
and $d_4 =\tfrac14 (d^4 -d^2 )$. Therefore $dd_3 -d_4 =\tfrac{1}{12}(d^4
-d^2 )$.

Now $\At{2}$ is generated by trees

\medskip
\centerline{
\hair=1pt
\labellist\small
\pinlabel {$a$} [br] at 19 76
\pinlabel {$b$} [tr] at 15 10
\pinlabel {$c$} [bl] at 159 76
\pinlabel {$d$} [tl] at 156 10
\endlabellist
\includegraphics[width=1.5in]{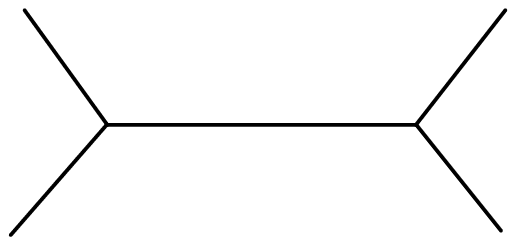}
}

\ni which we denote $T(a,b,c,d)$, where $a,b,c,d$ are elements of a basis $\B$ of $H$. Choose an ordering of $\B$.
\begin{lemma}\lbl{lem.at2}
$\At{2}$ is generated by $\{ T(a,b,c,d)\}$ with $a>b, c>d, a\ge c\ge b$ and, if $a=c$, $b\ge d$.
\end{lemma}

Assuming the lemma we can count the number of $T(a,b,c,d)$ satisfying the conditions of the lemma:
\begin{align*}
\binom{d}{4}&\qquad\text{ for }a>c>b>d  &
\binom{d}{4}&\qquad\text{ for }a>c>d>b  \\
\binom{d}{3}&\qquad\text{ for }a=c>b>d  &
\binom{d}{3}&\qquad\text{ for }a>c=b>d  \\
\binom{d}{3}&\qquad\text{ for }a>c>b=d  &
\binom{d}{2}&\qquad\text{ for }a=c>b=d 
\end{align*}
The sum of these six cases is exactly $\tfrac{1}{12}(d^4 -d^2 )$, which proves the Theorem.

\begin{proof}[Proof of \fullref{lem.at2}]
We first list some equalities:
\begin{align}
T(a,b,c,d)&= T(d,c,b,a) \lbl{eq.1}\\
T(a,b,c,d)&=- T(b,a,c,d)=-T(a,b,d,c)=T(b,a,d,c)\lbl{eq.2}\\
T(a,b,c,d)&=T(a,c,b,d)-T(a,d,b,c)\lbl{eq.3}\\
T(a,a,c,d)&=T(a,b,c,c)=0\lbl{eq.4}
\end{align}
Equation \eqref{eq.1} follows by rotating the tree. Equation \eqref{eq.2} is anti-symmetry and equation \eqref{eq.3} is the IHX relation. Equation \eqref{eq.4} follows from IHX (equation \eqref{eq.3}) for $c=d$, and then from \eqref{eq.1} for $a=b$.

We now prove the Lemma.  

It follows from equations \eqref{eq.2} and \eqref{eq.4} that $\At{2}$ is generated by $T(a,b,c,d)$ with $a>b$ and $c>d$. Next let $\T$ consist of all $\{ T(a,b,c,d)\}$ satisfying $a>b, c>d, a\ge c$ -- we show that $\T$ generates $\At{2}$. 

Suppose $a>b,c>d$ but $c>a$. Then we can use equations \eqref{eq.1} and \eqref{eq.2} to write $T(a,b,c,d)=T(c,d,a,b)$, which belong to $\T$. 

Now define $\T'$ to consist of all $\{ T(a,b,c,d)\}$ satisfying $a>b, c>d, a\ge c\ge b$. Suppose $T(a,b,c,d)\in\T$ but $b>c$. Then apply equation \eqref{eq.3} to $T(a,b,c,d)$ and note that $T(a,c,b,d)$ and $T(a,d,b,c)$ both belong to $\T'$.

It remains to eliminate those $T(a,b,c,d)\in\T'$ for which $a=c$ and $b< d$. Applying equations \eqref{eq.2}, \eqref{eq.3} and \eqref{eq.4} we have
\begin{align*}
T(a,b,a,d)&= T(a,a,b,d)-T(a,d,b,a) \\
 &= T(a,d,a,b)
\end{align*}  

This completes the proof.
\end{proof}

\bibliographystyle{gtart}
\bibliography{link}

\let\newaddresses\empty
\def\addresses{\address{}{Department of Mathematics\\
Brandeis University\\\newline
Waltham MA 02454-9110\\
USA\\\newline
\vspace{0pt}\\\newline
{\rm Correspondence to:\qua K Orr\\\newline}
Department of Mathematics\\Indiana University\\\newline
831 East 3rd St\\Bloomington IN 47405-7106\\USA
}}

\end{document}